\documentclass[12pt]{article}
\usepackage{a4,amssymb,theorem}
\ifx\mathbb\undefined\def\mathbb{\Bbb}\fi
\newtheorem{Theorem}{Theorem}

\newtheorem{Lemma}[Theorem]{Lemma}
\newtheorem{Proposition}[Theorem]{Proposition}
\newtheorem{Corollary}[Theorem]{Corollary}
\theorembodyfont{\normalfont}
\newtheorem{Proof}{Proof}

\newtheorem{Question}[Theorem]{Question}
\newtheorem{Definition}{Definition}
\makeatletter
\def\frak{\mathfrak}
\def\im{\mathop{\operator@font im}\nolimits}
\def\red{{\mathop{\operator@font red}\nolimits}}
\def\arg{\mathop{\operator@font arg}\nolimits}

\def\pole{\mathop{\operator@font pole}\nolimits}
\def\tor{{\mathop{\operator@font tor}\nolimits}}
\def\Cst{\mathop{\operator@font Cst}\nolimits}
\def\Newton{\mathop{\operator@font Newton}\nolimits}
\def\Ext{\mathop{\operator@font Ext}\nolimits}
\def\End{\mathop{\operator@font End}\nolimits}
\def\Ass{\mathop{\operator@font Ass}\nolimits}
\def\char{\mathop{\operator@font char}\nolimits}
\def\rk{\mathop{\operator@font rk}\nolimits}
\def\ker{\mathop{\operator@font ker}\nolimits}
\def\coker{\mathop{\operator@font coker}\nolimits}

\def\m{\mathop{\operator@font m}\nolimits}

\newcommand{\htt}{\mathop{\operator@font ht}\nolimits}
\newcommand{\ext}{\mathop{\operator@font ext}\nolimits}
\newcommand{\Supp}{\mathop{\operator@font Supp}\nolimits} 
\newcommand{\supp}{\mathop{\operator@font supp}\nolimits}
\newcommand{\Spec}{\mathop{\operator@font Spec}\nolimits}
\newcommand{\Edp}{\mathop{\operator@font E\mbox{-}dp}\nolimits}
\newcommand{\prof}{\mathop{\operator@font dp}\nolimits}
\newcommand{\Hom}{\mathop{\operator@font Hom}\nolimits}
\newcommand{\gr}{\mathop{\operator@font gr}\nolimits}
\newcommand{\h}{{\operator@font h}}
\newcommand{\LH}{{\operator@font H}}
\def\Ann{\mathop{\operator@font Ann}\nolimits}
\def\gr{\mathop{\operator@font gr}\nolimits}
\def\Tor{\mathop{\operator@font Tor}\nolimits}

\makeatother
\begin{document}
\sloppy
\title{Stiffness of Finite Free Resolutions and the Canonical Element Conjecture}
\author{Anne-Marie Simon and Jan R. Strooker}
\maketitle

\begin{abstract}
Sur un anneau local noeth\'{e}rien, les r\'{e}solutions libres minimales finies
ont une certaine propri\'{e}t\'{e} de rigidit\'{e} lorsque la 
caract\'{e}ristique de l'anneau \'{e}gale la caract\'{e}ristique de son corps 
r\'{e}siduel (th\'{e}or\`{e}me 1). 
Cette propri\'{e}t\'{e} \'{e}voque le 
crit\`{e}re de Buchsbaum et Eisenbud o\`{u} il n'est pas question de 
caract\'{e}ristique; cependant, personne ne sait si rigidit\'{e} reste valable 
en caract\'{e}ristique mixte. Conjecturons seulement ceci: sur tout anneau local
de Gorenstein, d'\'{e}gale caract\'{e}ristique ou non, toutes les r\'{e}solutions libres minimales finies sont rigides (au sens de la d\'{e}finition 1).

L'objet de cette note est de montrer, au th\'{e}or\`{e}me 2, que la conjecture pr\'{e}c\'{e}dente, qui est en fait une conjecture d'alg\`{e}bre lin\'{e}aire
pour une classe restreinte d'anneaux, \'{e}quivaut a la conjecture de l'\'{e}l\'{e}ment canonique de Hochster qui, elle, porte sur toute la classe 
des anneaux locaux noeth\'{e}riens. Ceci est fait via l'approximation par les modules de Cohen-Macaulay maximaux et les enveloppes de dimension injective 
finie. Cette th\'{e}orie, due a Auslander et Buchweitz, s'applique notamment 
aux modules de type fini sur un anneau de Gorenstein.

MSC: 13D02, 13D22, 13H10, 13C14, 13D05, 13D07, 13D25.                             

Keywords: finite free resolutions; stiffness; syzygies of finite projective 
dimension; Canonical Element Conjecture; Gorenstein rings; Auslander's $\delta$-invariant; hulls of finite injective dimension.
\end{abstract}

\section{Introduction}
	Throughout, $(A,\frak m, k)$ stands for a (commutative) noetherian
local ring with its maximal ideal and residue class field. We'll only deal with
such rings, and all modules will be finitely generated unless otherwise stated. For a $d$-dimensional ring, let $\bf F$ be a free resolution of $k$. Let $\bf x$ = $x_{1}, ..., x_{d}$ be any system of parameters in $A$; the surjection of $A/(x_{1}, ..., x_{d})$ onto $k$ lifts to a map from the Koszul complex ${\bf K}({\bf x},A)$ to $\bf F$. One says that the ring $A$ satisfies the Canonical Element Conjecture (CEC) provided the map $K_{d}({\bf x},A) \rightarrow F_{d}$ is never null. 
	
This conjecture was introduced by M. Hochster in [11] and proved for
equicharacteristic rings. Thanks to Heitmann's recent breakthrough, the result
is known in the unequal characteristic case when dim $A\leq 3$ [9]; also various special cases have been proved. Furthermore, there are several versions of CEC, which boil down to the same thing. Finally, Hochster and others have shown that CEC is equivalent with a number of conjectures, which at first sight appear rather distinct. The present note will not take up any of this; instead, it will add one more item to the list.

In order to prepare for the other notion in the title, consider an $A$-homomorphism $f: A^{m}$$\rightarrow A^{n}$ between free modules. Choosing bases for these, we can describe $f$ as an $n \times m$ matrix with coefficients in $A$. As an example, take $A = k[[X,Y]]/(XY)$ with $k$ a field of characteristic $\neq 2$. On certain bases, suppose the map $f: A^{2}$$\rightarrow A$ is given by the matrix $(x,y)$, where the lower case letters denote the images of $X$ and $Y$ in $A$. A change of base in $A^{2}$ achieves the row matrix $(x+y,x-y)$ as  description of $f$. Both $x$ and $y$ are zero divisors in $A$, but $x+y$ and $x-y$ are not.

In such a rectangular matrix, the elements in each particular column generate an ideal $\frak c$ $\subset A$ which we call a column ideal. We write $\gr \frak c$ for the length of a longest regular sequence contained in $\frak c$; in case $\frak c = A$, we put $\gr \frak c = \infty$. In our example, the grade depends on base choice. This motivates     
\begin{Definition} Let 
\begin{displaymath}
{\bf F} \quad = \quad 0\rightarrow F_{s} \rightarrow \cdots \rightarrow F_{i} \stackrel{d_{i}}{\rightarrow} F_{i-1} \rightarrow \cdots \rightarrow F_{0} 
\end{displaymath}be an acyclic complex of free's over $(A,\frak m,k)$
which is minimal in the sense that $k\otimes_{A}\bf F$ has null differentials. 
We say that $\bf F$ is `stiff' if, regardless of base choice, $\gr \frak c
\geq i$ for every column ideal belonging to $d_{i}$ for $i = 1, ..., s$.
We call the ring $A$ stiff if every such complex over it is stiff.
\end{Definition}

We prove
\begin{Theorem}\label{Theorem 1} Every noetherian local ring of equal characteristic is stiff.
\end{Theorem}
\begin{Theorem} Consider rings $(A,\frak m,k)$ of a fixed residual characteristic. Then every ring of this type satisfies \textup{CEC}
if and only if every Gorenstein ring of this type is stiff.
\end{Theorem}

Half a year after we had done this work, we belatedly discovered that 
Theorem~\ref{Theorem 1}, albeit in a different formulation, had already been obtained in [12, Th. 6.12] and [7, Th. 2.4]. The first of these articles uses the Frobenius endomorphism in finite characteristic and `lifts' to characteristic 0 using Artin approximation. This argument is sketched in the
second paper, but then another proof is given using big Cohen-Macaulay modules. 
Our proof in essence was the same as the latter. However, we believe it is 
more direct and transparent, and include it for this reason. 
\section{Preliminaries}
Before embarking on the proof of Theorem 1, we need to establish a few notational conventions and preliminary facts. Many of these topics were discussed in [17], so we use this as a blanket reference. 

Let $(A,\frak m, k)$ be our usual ring, which we drop from notation unless confusing. Since we shall also deal with infinitely generated modules, it is convenient to work with Ext-depth. Write $\ext^{-}(X,Y) = \inf\{n: \Ext^{n}(X,Y) \neq 0\}$ for two modules $X$ and $Y$; this number is a
nonnegative integer or $\infty$, as in other cases where we follow the same convention. For any $A$-module $X$ we write $\Edp X = \ext^{-}(k,X)$; when $X$ is finitely generated, this is equal to the usual depth of $X$ by way of regular sequences, for which we now write $\prof X$. Let $\frak p$ be a prime ideal in $A$, and put $k(\frak p)$ for the residue field of the local ring $A_{\frak p}$. Then there are the Bass numbers $\mu^{n}(\frak p,X) = \dim_{k(\frak p)}\Ext^{n}_{{A}_{\frak p}}(k(\frak p),X_{\frak p})$ [17, Th. 3.3.3]. Also
$\mu^{n}(\frak p,X) = \mu^{n}(\frak pA_{\frak p},X_{\frak p})$. 
Bringing into play local cohomology $\LH ^{n}_{\frak m}$, we find $\h ^{-}_{\frak m}(X) = \Edp X = \mu^{-}(\frak m,X)$, the first identity being taken care of by [17, Prop. 5.3.15]. If $\Edp X < \infty$, then $\Edp X\leq \dim A$ [17, Cor. 10.2.2], and if $A$ maps onto any ring $B$, then $\Edp_{A} Y = \Edp_{B} Y$ for every $B$-module $Y$ [17, Prop. 4.1.4]. We also write 
$\tor_{+}(X,Y) = \sup\{n: \Tor_{n}(X,Y) \neq 0\}$; this time it is a nonnegative integer or $-\infty$.

	We record an immediate consequence of the Acyclicity Lemma, e.g. [17, Th. 6.1.1], as  
\begin{Lemma} In a non-exact complex of modules, not necessarily finitely 
generated,  
$\bf X$ = $0 \rightarrow X_{m} \rightarrow \cdots \rightarrow X_{1} \rightarrow X_{0}$, for every $i \geq 1$ let either $\LH_{i}({\bf X}) = 0$ or $\Edp \LH_{i}({\bf X}) = 0$, and also  $\Edp X_{i} = q$. Then $q<m$.  
\end{Lemma}

The next inequality has presumably been noticed before. 
\begin{Lemma} For any of our rings $A$, let $\frak p$ be a prime ideal.
Then  $$\prof A \leq \prof A_{\frak p} + \dim A/\frak p.$$    
\end{Lemma}
\begin{Proof} Put $r$ = dp $A_{\frak p}$, $t$ = dim $A/\frak p$. Let 
$\frak p$ = $\frak p_{0}$ $\subset ... \subset \frak p_{t} = \frak m$ be a
strict saturated chain of prime ideals. Then $\mu^{r}(\frak p,A) = \mu^{r}(\frak pA_{\frak p},A_{\frak p}) \neq 0$ and iterated use of Bass' Lemma [17, Prop. 3.3.4] gives us $\mu^{r+t}(\frak m,A) \neq 0$. Therefore \\dp $A \leq r+t$.
\end{Proof}

A big Cohen-Macaulay module $C$ over $A$ is defined to be a not necessarily 
finitely generated module on which some system of parameters in $A$ forms a
regular sequence. Then $C \neq \frak m C$ and E-dp $C$ = dim $A$ for such a
module. A big Cohen-Macaulay module $C$ is called `balanced' if every system of parameters in $A$ is regular on $C$. Hochster constructed big Cohen-Macaulay 
modules whenever $A$ is equicharacteristic [10, Th. 5.1], and also constructed balanced ones. Bartijn and Strooker showed that the completion of
a big Cohen-Macaulay module with respect to the maximal ideal of $A$ is always
balanced [2, Th. 1.7 \& Th. 2.6], [17, Th. 5.2.3].

Dealing with infinitely generated modules, we need the small support, a selection of primes at which localization behaves decently. Let $X$ be an $A$-module, then $\supp X = \{\frak p \in \Spec A: \Edp_{A_{\frak p}} X_{\frak p} < \infty\}$. Small supports go back to Foxby [8, \S 2] who used $\tor_{-}$ rather than $\ext^{-}$, but this amounts to the same thing [17, Cor. 6.1.10]. Sharp showed in [14, Th. 3.2] that, for a balanced big Cohen-Macaulay module, $\supp$ has good properties with respect to height and dimension. For our purposes we record

\begin{Proposition} Let $C$ be a balanced big Cohen-Macaulay module over $A$
and let $\frak p \in \supp C$. Then
	$$\Edp_{A_{\frak p}}C_{\frak p} = \htt \frak p = \dim A_{\frak p} 
	\quad\mbox{ and }\quad \htt \frak p + \dim A/\frak p = \dim A,$$ while 
	$\frak p C_{\frak p} \neq C_{\frak p}.$
\end{Proposition} 

\section{Towards Theorem 1} First we sharpen a result of Evans-Griffith [6, Th. 1.11], replacing their $\htt \frak p$ by $\prof A_{\frak p}$. It seems also that the original proof required some mild condition on the ring.
\begin{Theorem} Let $M$ be a finitely generated $A$-module of finite projective 
dimension and $\frak p$ a prime ideal in $A$. Assume that C is a balanced big
Cohen-Macaulay module over $A/\frak p$. Then $\tor^{A}_{+}(C,M) \leq \prof A_{\frak p}$.
\end{Theorem} 
\begin{Proof}Since $C \neq \frak m C$ and $C\otimes_{A}M \neq 0$, putting 
$t = \tor_{+}(C,M)$, we see that $0 \leq t \leq \prof A$, as $t$ does not
exceed the projective dimension of $M$ which in turn is not greater than 
$\prof A$ by the Auslander-Buchsbaum identity. Let $\frak q$ be a prime which is minimal in $\Supp \Tor_{t}(C,M)$, so that $\frak q \supset \frak p$. Thus the 
nonnull $A_{\frak q}$-module $\Tor^{A_{\frak q}}_{t}(C_{\frak q},M_{\frak q}) =
\Tor^{A}_{t}(C,M)_{\frak q}$ is only supported in the maximal ideal of the 
local ring $A_{\frak q}$, so $\Edp_{A_{\frak q}}$ of this module is 0. 

Let
\begin{displaymath}
{\bf F}\quad =\quad 0\rightarrow F_{s}\rightarrow \cdots \rightarrow F_{1} \stackrel{ d_{1}}{\rightarrow} F_{0}
\end{displaymath}
be a minimal free $A_{\frak q}$-resolution of the module $M_{\frak q} = \coker
d_{1}$. This module's projective dimension is $s \leq \prof A_{\frak q}$.
By our choice of $\frak q$, we also know that $t\leq s$.

Now consider the complex $C_{\frak q}\otimes_{A_{\frak q}}{\bf F}_{\frak q}$ which we truncate at degree $t-1$, writing $F_{-1} = 0$ in case $t=0$. Its homology is concentrated in degree $t$, where its Ext-depth is 0. The Ext-depth of all its chain modules is $\Edp_{A_{\frak q}}C_{\frak q}$, so Lemma 3 tells us that $\Edp_{A_{\frak q}}C_{\frak q} \leq s-t < \infty$. In particular, $\frak q \in \supp C$.

We apply Proposition 5 to the ring $A/\frak p$ to obtain $\Edp_{A_{\frak q}}C_{\frak q} = \Edp _{(A/\frak p)_{\frak q}}C_{\frak q} = \dim (A/\frak p)_{\frak q} = \dim A_{\frak q}/\frak p A_{\frak q}$. Next we invoke 
Lemma 4 for the ring $A_{\frak q}$ and find $\prof A_{\frak q} \leq \prof A_{\frak p} + \dim A_{\frak q}/\frak p A_{\frak q}$. 

Tying all this together we establish the result
\begin{displaymath}
t \leq s - \Edp_{A_{\frak q}}C_{\frak q} = s - \dim A_{\frak q}/\frak p A_{\frak q} \leq \prof A_{\frak q} - \dim A_{\frak q}/\frak p A_{\frak q} \leq \prof A_{\frak p}.
\end{displaymath}
\end{Proof}
	
The statement we want is now an easy corollary.

\begin{Proof}[of Theorem 1] Let ${\bf F}$ be a complex as in Definition 1 over an equicharacteristic ring $A$, and consider its $i$-th boundary map $d_{i}$,
$1 \leq i \leq s$. Let $\frak c$ be, say, the first column ideal in some matrix description of $d_{i}$. Take a prime ideal $\frak p \supset \frak c$ with 
$\gr \frak c = \prof A_{\frak p}$. There exists a balanced big Cohen-Macaulay module $C$ over $A/\frak p$. The $i$-th boundary map of the complex 
$C\otimes_{A}{\bf F}$ maps the first copy of $C$ in $C\otimes_{A}F_{i}$ to 0
since $\frak c C$ = 0. This first copy itself is not in the image of the previous map because the complex $\bf F$ is minimal and $C\neq \frak m C$. Putting $M = \coker d_{1}$, we see that $\Tor^{A}_{i}(C,M) \neq 0$. With 
Theorem 6 then $i \leq \prof A_{\frak p} = \gr \frak c$.
\end{Proof}

At this point, what can one do in mixed characteristic? If $p$ is the 
residual characteristic, we can replace a column ideal $\frak c$ belonging
to the $i$-th map $d_{i}$ by the ideal $\frak c + (p) = \frak a$ and take a prime ideal $\frak q$ containing it such that $\prof A_{\frak q} = \gr \frak a$.
Tensoring $\bf F$ with a balanced big Cohen-Macaulay module over $A/\frak q$, and taking into account that $\gr \frak c \leq \gr \frak a \leq \gr \frak c +1$,
we obtain 
\begin{Proposition} Let $\bf F$ be a complex as in Definition 1. Then the 
column ideals belonging to $d_{i}$ always have grade at least $i-1$.
\end{Proposition}
\section{Stiffness depends on first syzygies}
In order to prepare for Theorem 2, we explain what we mean by this 
title and provide a straightforward proof of the appropriate statement.

\begin{Definition} Let $\Gamma$ be a class of rings such that, if $A\in \Gamma$
and $x\in \frak m$ is a non zerodivisor, then $A/(x)\in \Gamma$. We call such a 
class `consistent'.
\end{Definition}

Examples of consistent classes would be all rings of positive residual 
characteristic, the class of all Gorenstein rings, or their intersection.
\begin{Proposition} Let $\Gamma$ be a consistent class of rings. Then every
$A \in \Gamma$ is stiff if and only if for every complex $\bf F$ as in Definition 1 over every such ring, $\gr \frak c \geq 1$ for every column ideal 
belonging to the first boundary map $d_{1}$, no matter what bases.
\end{Proposition}
\begin{Proof} Assume the second condition and let $\bf F$ be a complex of length $s$ as in Definition 1 over a ring $A$ in $\Gamma$. We have to show that $\bf F$ is stiff.

For $s=1$ this is our assumption. So take $s\geq2$ and suppose that the result has been proved for complexes of lesser length. Let $\frak c$ be a column ideal in some matrix description of $d_{i}$, $2\leq i\leq s$. By our assumption, this 
ideal contains a non zerodivisor $x\in \frak m$. Write $M = \coker d_{1}$; then
$\Tor^{A}_{j}(A/(x),M) = 0$ for $j\geq 2$. Therefore, when we tensor with $A/(x)$ and clap bars on the $F_{i}$'s and $d_{i}$'s, 
\begin{displaymath}0 \rightarrow  \bar{F}_{s} \rightarrow \cdots \rightarrow \bar{F}_{2} \stackrel{\bar{d}_{2}}\rightarrow \bar{F}_{1}
\end{displaymath}
becomes a minimal acyclic complex of free's over the ring $A/(x)$. The ideal 
$\frak c/(x)$ is a column ideal of the ($i-1$)-st map in this new complex. The 
induction hypothesis allows us to conclude that its grade is at least $i-1$. 
But $x$ was chosen in $\frak c$, so $\gr \frak c \geq i$, and we are done.
\end{Proof}  

This proposition tells us that  `in bulk' the column ideals of the first boundary map are in control of stiffness. A form of `reduction to first syzygies' is already present in [7, Lemma 2.5].

In the proposition, there is a sufficient condition that $\gr \frak c \geq 1$ for each column ideal $\frak c$ in any matrix description of $d_{1}$ for every complex $\bf F$ as in Definition 1. For use in the next section, we state this as: in any syzygy $Z$ of finite projective dimension, $\Ann z = 0$ for every minimal generator $z$ of $Z$. Indeed, the columns of the matrix describe minimal generators of $Z = \im d_{1}$ as a syzygy in $F_{0}$, and the equivalence becomes evident.

\section{Towards Theorem 2}
In [16], we developed various concepts stemming from the seminal paper [1]
of M. Auslander and R.-O. Buchweitz, and showed how they intertwine with several of the Homological Conjectures, including CEC. Auslander-Buchweitz theory
essentially works for Cohen-Macaulay rings with dualizing module, and becomes even nicer for Gorenstein rings. For some time now, our point of view has been that these homological conjectures for all rings, can be interpreted as statements about modules over Gorenstein rings [18], [15], [16]. Thus `representation theory' of Gorenstein rings becomes pivotal. In this vein we quote from [15, \S 4] and [16, \S 6.4]
\begin{Theorem} The Canonical Element Conjecture is true for all rings of a 
certain residual characteristic, if and only if for every unmixed nonnull ideal $\frak b$ of zerodivisors in a Gorenstein ring $R$ of this residual characteristic, one has $\delta(R/\frak b) = 0$.
\end{Theorem}

Here $\delta$ is Auslander's invariant. For a module $M$ over a Gorenstein ring
$R$, the expression $\delta (M) = 0$ means that there exists a maximal 
Cohen-Macaulay module over $R$, which has no free direct summand and which 
surjects onto $M$. In general, $\delta (M) = n$ means that among all maximal 
Cohen-Macaulay modules mapping onto $M$, the integer $n$ is the smallest rank 
of their free summand which occurs, e.g. [16, Prop. 4.2]. In the proof of the above theorem, however, a new and somehow related characterization of the $\delta$-invariant is crucial. This time look at all possible surjections 
$p: R^{t} \rightarrow M$ from finitely generated free modules onto $M$. Put
$d$ for the dimension of $R$, then $\Ext^{d}_{R}(k,R^{t})$ is a $t$-dimensional
vector space over $k$. For every such $p$, the dimension of its image under
$\Ext^{d}_{R}(k,p)$ equals $\delta(M)$ [16, Th. 4.1 (iv)]. Covertly, this 
description of $\delta$ also plays a role in our treatment of the next result.

Before stating and proving this, we recall a few facts about modules
over a Gorenstein local ring $R$ which are just about standard. A module $M$ is a syzygy iff it is torsionless iff $\Ass M \subset \Ass R$. For a cyclic module 
$R/\frak b$ this is the case iff $\frak b$ consists of zero divisors and is
unmixed iff $\frak b$ is an annihilator ideal, i.e. $\frak b = \Ann \frak a$
for some ideal $\frak a \subset R$. We also indiscriminately speak of finite
injective dimension or finite projective dimension, since over a Gorenstein ring
these designations apply to the same modules, e.g. [17, Th. 10.1.9]. This allows for a quick proof of
\begin{Proposition} Let $R$ be a Gorenstein ring. Then $\delta (R/\frak b) = 0$
for every unmixed nonnull ideal of zero divisors $\frak b$ if and only if 
$\Ann z = 0$ for every minimal generator $z$ of every syzygy of finite projective dimension.
\end{Proposition}   
\begin{Proof} Suppose the condition on the $\delta$'s is satisfied, and let
$Z$ be a syzygy of finite projective dimension and $z \in Z$. 
Since $z$ also lives in a free module, $\Ann z = \frak b$ is an annihilator ideal, so that $\delta(Rz) = 0$ unless $\frak b = 0$. If not, the injection
of $Rz$ into $Z$ takes $z$ into $\frak m Z$ by [16, Prop. 4.5 (iii)], so 
$z$ is not a minimal generator.   

Now assume the condition on minimal generators, and let $\frak b \subset R$
be a nonnull unmixed ideal of zero divisors. Take a hull of finite injective dimension of
$R/\frak b$ [1,\S 0], [16, \S 1 \& \S 3]. In other words, there is a short exact sequence
\begin{displaymath}
		0 \rightarrow R/\frak b \rightarrow Z \rightarrow C \rightarrow 0
\end{displaymath}
where $Z$ has finite injective dimension and $C$ is maximal Cohen-Macaulay. 
Since $\Ass Z \subset \Ass R/\frak b \cup \Ass C$ and all the primes on the 
right hand side belong to $\Ass R$, we see that $Z$ is a syzygy. By assumption, $R/\frak b$ cannot hit a minimal generator of $Z$, so lands in $\frak m Z$, which means that  $Z/\frak m Z \simeq C/\frak m C$. Hence $\delta(R/\frak b) = 0$ by [16, Prop. 4.4].
\end{Proof}

We are poised to clinch the second claim in the Introduction.

\begin{Proof}[of Theorem 2] Combine Theorem 9, Propositions 8 and 10  with the remark at the end of section 4 which interprets columns as minimal generators of syzygies. 
\end{Proof}

With Theorem 1 and the fact that CEC holds in equal characteristic, one notices that the word `residual' in the statement of Theorem 2 may be replaced by `mixed'. 
\section{More about stiffness}
Very recently, we became aware that Hochster and Huneke had continued their 
investigations on topics around stiffness after [12]. In their long and 
exciting article [13], \S 10.7-10.12 contain extensions of their earlier results  and of [7] which we briefly present from our point of view. 

Their work, like that of Evans and Griffith, is based on order ideals. Let
$M$ be an $A$-module, and $m\in M$. Put $M^{*} = \Hom_{A}(M,A)$. Then $M^{*}(m)$ is the ideal which is generated in $A$ by the images of $m$ under all the maps
in $M^{*}$. This order ideal used to be denoted by $O_{M}(m)$, but fashion appears to be changing. Their main result in this direction [13, Th. 10.8]
reads: let $A$ be of positive characteristic, and let $z$ be a minimal 
generator in an $i$-th syzygy $Z$ of finite projective dimension. Then 
$\gr Z^{*}(z) \geq i$. We argue that every stiff ring has this property.

Indeed, let $\bf G$ be a finite free resolution of some $M$ in which $Z$ is an $i$-th syzygy of which $z$ is a minimal generator. By an old but basic result of Eilenberg [5, Th. 8], $\bf G$ is a direct sum of a minimal free resolution 
$\bf F$ of $M$ and a free resolution $\bf H$ of 0. Then $Z = U \oplus V$ where 
$U$ and $V$ are $i$-th syzygies in $\bf F$ resp. $\bf H$. Let $z=u+v$ be the corresponding decomposition. Since $\frak m Z = \frak m U \oplus \frak m V$, 
$u$ and/or $v$ needs to be a minimal generator of the syzygy to which it belongs. In case $V \neq 0$ and $v$ is a minimal generator, $\gr V^{*}(v) = \infty$ because $\bf H$ splits and $V$ is free. If $u$ is a minimal generator of $U$, then by stiffness its column ideal $\frak c$ has grade $\geq i$ for any choice of bases in the complex $\bf F$. But $\frak c = F_{i-1}^{*}(u) \subset U^{*}(u)$. By the canonical projections of $Z$ onto its components one sees that 
$Z^{*}(z) \supset U^{*}(u) + V^{*}(v)$, so $\gr Z^{*}(z) \geq i$ in all cases.

For the next two theorems we first derive a couple of immediate consequences of 
stiffness. In a stiff complex $\bf F$ of length $s$ take $i < s$. Then 
$\gr \frak c \geq i$ for all column ideals belonging to $d_{i}$, regardless of bases, so such an ideal contains a regular sequence of length at least $i$. 
The number of nonnull generators of the ideal is at most the number of rows in the matrix $X$ depicting $d_{i}$, say $f_{i-1} = \rk F_{i-1}$. Therefore 
$f_{i-1} \geq i$. Since also $f_{i} \geq i+1$, we see that $Z = \im d_{i}$
requires at least $i+1$ generators and also that $i \times i$ minors occur in $X$. In Theorem 14 we shall denote by $\frak c_{i,t}, 1 \leq t \leq i$, the ideal generated by all the $t \times t$ minors which occur in $t$ fixed columns of $X$. Notice that, and this is the point of the exercise, we are not concerned with bases nor with which $t$ columns we concentrate on.

In [13, Cor. 10.10] a considerable strengthening of the Evans-Griffith syzygy theorem [6, Cor. 3.16] is proved for rings of positive characteristic using [13, Th. 10.8]. With the above, one can state 
\begin{Theorem} Let $Z$ be an i-th syzygy of finite projective dimension over a stiff ring. If $Z$ is not free, then it requires at least $i+1$ generators.
In any case, if $z_{1}, ... , z_{t}$, $1 \leq t \leq i$, form part of a minimal system of generators for $Z$, then they generate a free submodule G of $Z$.
Moreover the factor module Z/G is an ($i-t$)-th syzygy.  
\end{Theorem}

The next result we aim for is a stiff version of [13, Cor. 10.11].
Since we'll offer a slightly different proof, we make some preparations. 
Among ways to recognize syzygies of finite projective dimension, this one 
is useful [13, Lemma 10.9]:
\begin{Lemma} For a module $M$ of finite projective dimension over a ring $A$,
the following conditions are equivalent.
\begin{enumerate}
\item[(i)] $M$ is an $i$-th syzygy;
\item[(ii)]	For every prime ideal $\frak p$ in $A$, either 
	$\prof M_{\frak p} \geq i$ or $M_{\frak p}$ is free.  
\end{enumerate}
\end{Lemma}
\begin{Corollary} If in a short exact sequence of modules the two outer terms are $i$-th syzygies of finite projective dimension, so is the middle term.
\end{Corollary}

Next recall the celebrated Buchsbaum-Eisenbud criterion [4, Th.], [17, Th. 6.2.3] which tells us when a finite free complex is exact. Let $\bf F$ be such a complex of length $s$ with boundary maps $d_{i}$. For any map $\phi$ between free modules, put $I_{u}(\phi)$ for the ideal of $A$ generated by the $u \times u$ minors in a matrix description of $\phi$; this does not depend on choice of bases. The rank of $\phi$ then is the largest $u$ for which $I_{u}(\phi) \neq 0$. Put $r_{i} = f_{i} - f_{i+1} + ... \pm f_{s}, 1 \leq i \leq s$. Then $\bf F$ is acyclic if and only if $\gr I_{r_{i}}(d_{i}) \geq i$ for all $i$. In this case, moreover, $\rk d_{i} = r_{i}$ [3, Th. 1]. Observe that here minimality nor characteristic are mentioned.

Now comes the assertion we're after, which paraphrases [13, Cor. 10.11]. We
keep notation. 
\begin{Theorem} Let $\bf F$ be a complex of length $s$ as in 
Definition 1 over a stiff ring $A$. Then $\gr \frak c_{i,t} \geq i-t+1$ for
$1 \leq t \leq i < s$. 
\end{Theorem}
\begin{Proof} A set of $t$ columns of a matrix describing $d_{i}$ corresponds
to a part $z_{1}, ..., z_{t}$ of a minimal set of generators of $Z_{i} = 
\im d_{i}$. The submatrix formed by these columns describes the injection $g$
of the module $G = Az_{1} + \cdots Az_{t}$, which is free by Theorem 11, into $F_{i-1}$. There is an exact sequence $0 \rightarrow Z_{i}/G \rightarrow F_{i-1}/G \rightarrow Z_{i-1} \rightarrow 0$. In this sequence, the modules $Z_{i}/G$ and $Z_{i-1}$ are ($i-t$)-th syzygies, the former by Theorem 11. In view of Corollary 13, so is $F_{i-1}/G$. Hence $g$ is the last boundary map in an acyclic complex of free's of length $i-t+1$. We conclude with the Buchsbaum-Eisenbud criterion.
\end{Proof}
\begin{Question} For $t=1$, the conclusion of the theorem just reasserts stiffness. As $t$ grows, larger and larger submatrices come into play, but in $t$ fixed columns at a time. And their determinants generate smaller and smaller ideals. What is the connexion with the Buchsbaum-Eisenbud criterion, which we playfully employed in our proof? Aficionado's of multilinear algebra are kindly invited to shed light on this.
\end{Question}

\

\

\begin{tabular}{ll}
Anne-Marie Simon			&	Jan R. Strooker\\
Service d'Algebre C.P. 211		&	Mathematisch Instituut\\ 
Universite Libre de Bruxelles\qquad\qquad\qquad	&	Universiteit Utrecht\\
Campus Plaine				&	Postbus 80010\\
Boulevard du Triomphe			&	3508 TA Utrecht,\\ 
B-1050 Bruxelles, Belgique		&	The Netherlands\\
e.mail: amsimon@ulb.ac.be		&	e.mail: strooker@math.uu.nl  	
\end{tabular}

\end{document}